\documentstyle[11pt]{article}
\normalbaselines

\newcommand{\Complex}{\mbox{\sf I\hskip -1.7mm\rm C}}
\newcommand{\Integer}{\mbox{\bf Z}}
\newcommand{\Rational}{\mbox{\bf Q}}
\newcommand{\Natural}{\mbox{\bf N}}
\def\inO #1{\mbox{\bf O\hskip -7pt\raise 1.5pt\hbox{\scriptsize #1}}}
\def\inSO #1{\mbox{\bf SO\hskip -7pt\raise 1.5pt\hbox{\scriptsize #1}}}

\begin{document}
\bibliographystyle{unsrt}

\vbox{\vspace{38mm}}
\begin{center}
{\Large \bf Generalized Clifford Algebras  and the Last Fermat Theorem
 \\[5mm]}

 A. K. Kwa\'sniewski$^{*,**,***}$ ; W. Bajguz$^{**}$ \\[3mm]
{\it
    $^*$    Inst. of Math. Technical University of Bia{\l}ystok,
        ul. Wiejska 45 A, room 133\\
    $^{**}$ Inst. of Phys. Warsaw University Campus Bia{\l}ystok,
        ul. Przytorowa 2 A,\\
                15--104 Bia{\l}ystok, Poland\\
    $^{***}$    e-mail: Kwandr @ cksr.ac.bialystok.pl\\
}

\end{center}

\begin{abstract}
One shows that the Last Fermat Theorem is equivalent to the statement that
all rational solutions $x^k+y^k=1$ of  equation ($k\geq 2$) are provided by
an orbit of rationally parametrized subgroup of a group preserving $k$--ubic
form. This very group naturally arrises in the generalized Clifford algebras
setting \cite{K1985}.
\end{abstract}

{\bf I.} The stroboscopic motion of the independent oscilatory
degree of freedom is given by iteration of the "classical map" matrix
\begin{equation}
L(\Delta )=\frac{1}{1+\Delta^2}\left(
\begin{array}{cc}
1-\Delta^2 & -2\Delta \\
2\Delta & 1-\Delta^2
\end{array}
\right)\ \ \Delta\in\bar{\Rational}=\Rational\cup\{\infty\}
\label{1}
\end{equation}
(see \cite{BK} and references therein).

$L\left(\Delta\right)$ of (\ref{1}) provides the rational
parametrization of the unit circle obtained via stereographic projection
composed with $\pi/2$--rotation represented by an
${\bf i}=\left(\begin{array}{cc}
0 & -1 \\ 1 & 0
\end{array}\right)$ imaginary unit matrix.

The set $\inSO{2}\left( 2;\Rational\right) = \{ L(\Delta );
  \:\Delta\in\bar{\Rational}\}$
is the well known group. Ocasionally it is a refine exercise to prove that
$L(\Delta_1 )L(\Delta_2 )=L(\Delta )$, where
$\Delta_1,\Delta_2\in\bar{\Rational}$
\begin{equation}
\Delta=\frac{\Delta_1+\Delta_2}{1-\Delta_1\Delta_2}=
   \frac{\left( 1+\Delta_1^2\right)\left( 1+\Delta_2^2\right) -
     \left( 1-\Delta_1^2\right)\left( 1-\Delta_2^2\right) +4\Delta_1\Delta_2}
     {2\left[\Delta_1\left( 1-\Delta_2^2\right) +
       \Delta_2\left( 1-\Delta_1^2\right)\right]}
  \label{2}
\end{equation}
with special cases such as $L(1)L(1)=L(\infty )=-{\bf 1}$ or
$L(\Delta )L(-\Delta )=L(0)={\bf 1}$ included (for the last one use
d'Hospital rule).

One is tempted to call the group
\[ \inO{2}\left( 2;\bar{\Rational}\right) =
\inSO{2}\left( 2;\bar{\Rational}\right)\cup\left(\begin{array}{cc}
1 & 0 \\ 0 & -1 \end{array}\right)\inSO{2}\left(2;\bar{\Rational}\right) \]
the 2--Fermat group as it preserves quadratic form
\begin{equation}
x^2+y^2=1\ \ \ \ x,y\in\Rational \label{3}
\end{equation}
and even more.

\underline{Observation}: $\inO{2}\left( 2;\bar{\Rational}\right)$ group
acts transitively on the set of all rational solutions of (\ref{3}).

Proof: For any two solutions
$\left(\begin{array}{c}x_0\\y_0\end{array}\right)$,
$\left(\begin{array}{c}x\\y\end{array}\right)$
one easily finds $A\in\inO{2}\left( 2;\bar{\Rational}\right)$ such that
$\left(\begin{array}{c}x\\y\end{array}\right) =
 A\left(\begin{array}{c}x_0\\y_0\end{array}\right)$.
For example: let $x\neq-x_0$ and let $y\neq -y_0$; then $A=L(\Delta )$,
\[\Delta = \frac{x_0y-xy_0}{x_0(x_0+x)+y_0(y_0+y)}\]
The shape of formula for $\Delta$ depends on the way one chooses to find it
out. One way is just straightforward calculation. The other is based on the
observation that for
$\left(\begin{array}{c}x_0\\y_0\end{array}\right)\equiv
  \left(\begin{array}{c}1\\0\end{array}\right)$,
$\Delta =\frac{y}{x+1}$. Hence for any
$\left(\begin{array}{c}x_0\\y_0\end{array}\right)$\ \&
$\left(\begin{array}{c}x\\y\end{array}\right)$ the corresponding $\Delta$
is being found due to the obvious identity
$L\left(\Delta\right)\equiv
  L\left(\frac{y}{x+1}\right) L\left( -\frac{y_0}{x_0+1}\right)$.
That way we arrive at the intriguing identity valid for all {\bf solutions}
of (\ref{3}) i.e.
\begin{equation}
\frac{x_0y-xy_0}{x_0(x_0+x)+y_0(y_0+y)}\equiv
  \frac{x_0y-xy_0+y-y_0}{x_0(x_0+x)+y_0(y_0+y)+x+x_0} \label{Exy}
\end{equation}

\underline{Conclusion}: It is enough to start with trivial solution
$\left(\begin{array}{c}x_0\\y_0\end{array}\right) =
 \left(\begin{array}{c}1\\0\end{array}\right)$ of (\ref{3}). All others are
obtained as elements of the corresponding orbit of
$\inO{2}\left( 2;\bar{\Rational}\right)$ i.e.  2--Fermat group.

\underline{Remark}: An iteration of $L\left(\Delta\right)$, i.e.
$L\left(\Delta\right)\rightarrow L^2\left(\Delta\right)\rightarrow ...
L^k\left(\Delta\right)\rightarrow ...$
provides us with stroboscopic motion in one oscilatory degree of freedom
which in view of (\ref{2}) is chaotic; it is in a sense -- "number theoretic"
-- chaotic. (For the relation to Fibonacci--like sequences -- see \cite{BK})

\medskip

{\bf II.} Consider now
\begin{eqnarray}
x^k+y^k=1\ \ k\geq 3,\ n\in\Natural \label{4} \\
\mbox{where}\ \ x,y\in\Complex . \nonumber
\end{eqnarray}
Denote by $\inO{k}\left( 2;\Complex\right)$ the group of all
\underline{linear} transformations preserving this $k$--ubic form \cite{K1985}
related to generalized Clifford algebras \cite{K1985}. Of course starting
from any -- say trivial solution
$\left(\begin{array}{c}1\\0\end{array}\right)$
of (\ref{4}), the orbit $\inO{k}\left( 2;\Complex\right)$ would provide us with
a family of other solutions. Starting with another, nontrivial solution
\[\left(\begin{array}{c}x \\ \sqrt[k]{1-x^k}\end{array}\right)
   \ \ \ x\neq1\ \ x\in\Complex^\ast = \Complex\setminus\{ 0\}\]
we get -- for each another $x$ (not belonging to the precedent orbit!) a new
orbit of solutions. Evidently the set of all complex solutions of
$x^k+y^k=1$ has the structure of the sum of disjoint orbits of
$\inO{k}\left( 2;\Complex\right)$. In this connection note that the relation
between two solutions belonging to different orbits must be nonlinear.

According to K. Morinaga and T. Nono \cite{MN1952,K1985}
\begin{equation}\inO{k}\left( n;\Complex\right)=\left\{
\omega^l\delta_{i,\sigma (j)};\ l\in\Integer_k,\:\sigma\in\mbox{\bf S}_n
\right\}\ \ \ k\geq3 \label{5} \end{equation}
where $\omega=\exp\left\{\frac{2\pi i}{k}\right\}$. Naturally
$|\inO{k}\left( n;\Complex\right) |=k^nn!$, hence every "$k$--Fermat group"
orbit of solutions of (\ref{4}) counts $2k^2$ elements.

One readily notices that the orbit
$\inO{k}\left( 2;\Complex\right)
\left(\begin{array}{c}1 \\ 0\end{array}\right)$
does not exhibit any nontrivial rational solution, as the $k$--Fermat group,
$k\geq 3$ i.e. $\inO{k}\left( n;\Complex\right)$ contains the only one
{\bf rationaly} parametrized subgroup, i.e. the matrix permutation subgroup
$\simeq\mbox{\bf S}_n$.

Thus we arrive at the

\underline{Conclusion}: The Last Fermat Theorem is equivalent to the
statement, that all avaiable rational solutions of
$x^k+y^k=1\ \ k\geq 2$ are provided by the orbit
$\inO{k}\left( 2;\bar{\Rational}\right)
\left(\begin{array}{c}1 \\ 0\end{array}\right)$;\ \
$\inO{k}\left( n;\Rational\right)\subset\inO{k}\left( n;\Complex\right)$.

One is evidently tempted to conjecture the "corresponding Last Fermat Theorem"
concerning $\inO{k}\left( n;\bar{\Rational}\right)\ \ n>2$ group. Hence
$n$--\underline{hypothesis}. Let $n\geq 2$, then
\[x_1^k+x_2^k+...+x_n^k=1\]
has no rational solutions for $k\geq 3$, except for trivial ones, i.e.
$x_s=0,\pm 1$, $s=1,...,n$.

This is however obviously {\bf false}, since for each $x_1$ -- natural and $k$ --
odd numbers it is easy to find natural $n$ and $x_2,...,x_n$ such that equation
is true. Anyhow quadratic forms for $n=2$ (appropiate to associate oscilations
with!) seem to be the only ones among $k$--ubic forms ($k\geq 2,\ n=2$) that
would provide us with nontrivial stroboscopic motion by group element
iteration as outlined in \cite{BK}.

\underline{Remark 1}: $k$--ubic forms of (1,1) signature as well as
corresponding generalized Clifford algebras are at hand \cite{K1985}, hence
the "2--hypothesis" equipped with (1,1) signature is easy to formulate;
namely:

Let $Q$ be a $k$--ubic form of (1,1) signature. Let $\vec{x}\in\Rational^2$;
then the all solutions of $Q(\vec{x})={\bf 1}$ are given by the orbit
\[\inO{k}\left( 1,1;\Rational\right)
\left(\begin{array}{c}1 \\ 0\end{array}\right) .\]
(This is of course equivalent to the (2,0) signature case)

\medskip

For the sake of examplification take $k=2,\ n=2$. Then
\[\inO{2}\left( 1,1;\bar{\Rational}\right) =
  \inSO{2}\left( 1,1;\bar{\Rational}\right)\cup
  \left(\begin{array}{cc}1 & 0\\ 0 & -1 \end{array}\right)
  \inSO{2}\left( 1,1;\bar{\Rational}\right)\]
where
\[\inSO{2}\left( 1,1;\bar{\Rational}\right)\equiv
  \{\tilde{L}(\Delta );\ \Delta\in\bar{\Rational}\};\ \
  \tilde{L}(\Delta )\equiv\frac{1}{1-\Delta^2}\left(\begin{array}{cc}
    1+\Delta^2 & 2\Delta \\ 2\Delta & 1+\Delta^2\end{array}\right).\]
It is then easy to see, that

\underline{Observation}: $\inO{2}\left( 1,1;\bar{\Rational}\right)$ group
acts transitively on the set of all rational solutions of $x^2-y^2=1$.

Proof: For any two solutions
$\left(\begin{array}{c}x_0\\y_0\end{array}\right)$,
$\left(\begin{array}{c}x\\y\end{array}\right)$
one easily finds $L(\Delta )\in\inO{2}\left( 1;1;\bar{\Rational}\right)$
such that
$\left(\begin{array}{c}x\\y\end{array}\right) =
 L(\Delta )\left(\begin{array}{c}x_0\\y_0\end{array}\right)$.
For example: let $x\neq-x_0$ and $y\neq -y_0$; then one has the following identity
\begin{equation}
\Delta = \frac{xy_0-yx_0}{x(x_0+x)+y(y_0+y)}\equiv
\frac{x_0y-xy_0+y-y_0}{x_0(x_0+x)-y_0(y_0+y)+x+x_0} \label{HDeltaxy}
\end{equation}
(\ref{HDeltaxy}) is analogous to (\ref{Exy}) i.e. it is valid on the set of
solutions of "hiperbolic" Fermat $n=2$ equation
\[x^2-y^2=1\ ;\ \ x,y\in\bar{\Rational}\]

The formula analogous to (\ref{2}) has the form:
\begin{equation}
\Delta=\frac{\Delta_1+\Delta_2}{1+\Delta_1\Delta_2}=
   \frac{\left( 1-\Delta_1^2\right)\left( 1-\Delta_2^2\right) -
     \left( 1+\Delta_1^2\right)\left( 1+\Delta_2^2\right) -4\Delta_1\Delta_2}
     {2\left[\Delta_1\left( 1+\Delta_2^2\right) +
       \Delta_2\left( 1+\Delta_1^2\right)\right]}
\label{HDelta}
\end{equation}
where
\[ L\left(\Delta\right)\equiv L\left(\Delta_1\right) L\left(\Delta_2\right)\ ;
\ \ \ L\left(\Delta_1\right) ,\ L\left(\Delta_2\right)\in
\inSO{2}\left( 1,1;\bar{\Rational}\right)\]
with special cases such as
$\tilde{L}(\Delta )\tilde{L}\left( -\frac{1}{\Delta}\right) =
  \tilde{L}(\infty )=-{\bf 1}$
or $\tilde{L}(\Delta )\tilde{L}(-\Delta )=\tilde{L}(0)$ included.

\underline{Remark 2}: We suggest relevance of hyperbolic functions of
$k$--th order \cite{KC1992} in relations between LFT and generalized Clifford
algebras (as used to linearize $k$--ubic forms in a Dirac way).


\begin{thebibliography}{99}

\bibitem{K1985}A. K. Kwa\'sniewski, J. Math. Phys., {\bf 26}, 9 (1985), 2234.
\bibitem{BK}W. Bajguz, A. K. Kwa\'sniewski, {\it On Quantum Mechanics and
 Fibonacci Sequences}, Advances in Applied Clifford Algebras, Vol. 4 (1994),
 73--88;
\bibitem{MN1952}K. Morinaga, T. Nono, J. Sci. Hiroshima Univ., Ser. A. Math.
 Phys. Chem., {\bf 16}, 13 (1952).
\bibitem{KC1992}A. K. Kwa\'sniewski, {\it On Hyperbolic and Elliptic Mappings
 and Quasi--Number Algebras}, Advances in Applied Clifford Algebras, Vol. 2,
 No 1 (1994), 107--144.

\end{thebibliography}
\end{document}